\documentclass[11pt]{article}

\usepackage[latin1]{inputenc}

\usepackage{amssymb,amsmath,amsthm}
\usepackage{array}
\usepackage[dvipsnames]{xcolor}
\usepackage{enumitem}
\usepackage[colorlinks]{hyperref}
\usepackage[margin=2cm]{geometry}
\usepackage{xypic}
\usepackage{needspace}

\newtheorem{theo}{Theorem}[section]
\newtheorem{coro}[theo]{Corollary}
\newtheorem{prop}[theo]{Proposition}
\newtheorem{rem}[theo]{Remark}
\newtheorem{lemma}[theo]{Lemma}

\newtheoremstyle{algth}
{3pt}
{3pt}
{}
{}
{\bfseries}
{}
{.5em}
{}

\theoremstyle{algth}
\newtheorem{algorithm}[theo]{Algorithm}

\newenvironment{alg}[2]
{\vspace{6pt}
	\needspace{\baselineskip}\noindent\hrulefill
	\begin{algorithm}\hfill		
		
		\vspace{-5pt}
		\noindent\hrulefill
		\vspace{3pt}
		\begin{description}
			\item[Input:] #1
			\item[Output:] #2 
			\item[Steps:]
		\end{description} 
		\begin{enumerate}[label=\arabic*.]
		}
		{ 
		\end{enumerate}
	\end{algorithm}
	\noindent\hrulefill\needspace{\baselineskip}
	\vspace{6pt}
}

\newcommand{\Q}{\mathbb Q}
\newcommand{\R}{\mathbb R}
\newcommand{\Z}{\mathbb Z}

\newcommand{\C}{\mathbb C}
\newcommand{\FF}{\mathbb F}

\newcommand{\HH}{\mathbf H}

\newcommand{\cB}{\mathcal{B}}
\newcommand{\cT}{\mathcal{T}}
\newcommand{\cO}{\mathcal{O}}

\newcommand{\sfM}{\mathsf{M}}
\newcommand{\sfA}{\mathsf{A}}

\newcommand{\zetaCC}[1]{z_{#1}}

\newcommand{\eps}{\varepsilon}
\newcommand{\on}[1]{\operatorname{#1}}

\renewcommand{\tilde}{\widetilde}

\DeclareMathOperator{\id}{{id}}
\DeclareMathOperator{\divv}{div}
\DeclareMathOperator{\Sp}{Sp}

\DeclareMathOperator{\im}{Im}

\begin{document}

\title{An inverse Jacobian algorithm for Picard curves}

\author{Joan-C. Lario and Anna Somoza (appendix by Christelle Vincent)}

\date{\today}

\maketitle 

\abstract{We study the inverse Jacobian problem for the case of Picard curves over $\C$. More precisely, we elaborate on an algorithm that, given a small period matrix $\Omega\in \mathbb{C}^{3\times 3}$ corresponding to a principally polarized abelian threefold equipped with an automorphism of order $3$, returns a Legendre--Rosenhain equation for a Picard curve with Jacobian isomorphic to the given abelian variety.
	
Our method corrects a formula obtained by Koike--Weng in \cite{KW} which is based on a theorem of Siegel. As a result, we apply the algorithm to obtain (numerically) all the isomorphism classes of Picard curves with maximal complex multiplication attached to the sextic CM-fields with class number at most \(4\). In particular, we obtain (conjecturally) the complete list of CM Picard curves defined over $\mathbb{Q}$.}

In the appendix, Vincent gives a correction to the generalization of Takase?s formula for the inverse Jacobian problem for hyperelliptic curves given in \cite{genus3hyper}.

\section{Introduction}

Let \(J\) be the map from the set \(\sfM_g\) of isomorphism classes of algebraic curves of genus~\(g\) defined over \(\C\) to the set \(\sfA_g\) of isomorphism classes of  complex principally polarized abelian varieties of dimension~\(g\). In this context, the \emph{inverse Jacobian problem} consists of identifying the preimage via \(J\) of the class of a given principally polarized abelian variety, if it exists. This is a classic result in the case of curves of genus~1, and has also been solved for curves of genus~2 \cite{Rosenhain, Thomae} and genus~3 \cite{genus3hyper,Guardia, KW, Takase, Weber, W}. Note that in all these cases, the map \(J\) is a bijection.

\bigskip
In this paper we present an inverse Jacobian algorithm for the family of Picard curves. This was initially done by Koike and Weng in \cite{KW}, but their exposition presents some gaps and mistakes that we fix here.

In Section \ref{sec:thomae-P} we give a formula to approximate the \(x\)-coordinates of the affine branch points of a Picard curve in terms of theta constants of its Jacobian, see Theorem~\ref{thm:thomae-P}. The given formula differs from the result in \cite{KW} by a third root of unity, see Remark~\ref{rmk:ThomaevsKW}.

In Section \ref{sec:algorithm-P} we first characterize the image under \(J\) of this family of curves, and then develop the algorithm that takes the Jacobian of a Picard curve \(C\) and returns a Legendre--Rosenhain equation for~\(C\), see Algorithm~\ref{alg:Picard}. The main step of the algorithm is applying the formula of Theorem~\ref{thm:thomae-P}, so we first identify the objects needed to apply said formula, mainly the Riemann constant and the images by the Abel-Jacobi map of the affine branch points. Our algorithm makes the process of identifying these points explicit in Theorem~\ref{thm:setDs-P}, see Remark~\ref{rmk:findDsvsKW} for a comparison with the approach of \cite{KW}.

Our correction of the algorithm allows us to  re-obtain the results of \cite{KW} and extend the list of maximal CM Picard curves, that is, Picard curves such that their Jacobians have endomorphism ring isomorphic to the maximal order of a sextic CM number field~\(K\). We obtain twenty-three new curves, displayed in Section~\ref{Llista}, among which we include all maximal CM Picard curves defined over \(\Q\). The corresponding CM-fields are collected from \cite{Klist}. The computations have been performed using SageMath~\cite{SageMath}, and an implementation can be found at \cite{github-algorithms}.

In the appendix, Vincent applies the tools introduced in Section~\ref{sec:thomae-P} to correct a sign in the generalization of Takase's formula for the inverse Jacobian problem for hyperelliptic curves, given in a previous article \cite{genus3hyper}.

The present paper is an extension and clarification of our earlier work \cite{LS16} to include further improvements of the algorithm, such as Theorem~\ref{thm:setDs-P}.

\subsection*{Acknowledgements}
The authors would like to thank Marco Streng and Christelle Vincent for useful discussions.

\section{A Thomae-like formula for Picard curves\label{sec:thomae-P}}

Let \(C\) be a Picard curve defined over \(\C\), that is, a genus-3 smooth, plane, projective curve given by the affine equation \(y^3 = f(x)\) where \(f\) is a polynomial of degree~4. The curve \(C\) has an automorphism~\(\rho\) of order~3 given by \((x,y)\mapsto (x,\zetaCC{3}y)  \) with \(z_3 = \exp\left(\frac{2\pi i}{3}\right)\). This automorphism fixes the \emph{affine branch points} \((t,0)\) with \(f(t) = 0\). The curve \(C\) has a unique point at infinity, with projective coordinates \((0:1:0)\), which is also fixed by the automorphism \(\rho\).

Up to isomorphism, we can (and do) assume that \(C\) is given by a \emph{Legendre--Rosenhain equation} 
\begin{equation}\label{eq:LR-Picard}
y^3 = x(x-1)(x-\lambda)(x-\mu).
\end{equation}

Following the literature, for example \cite[Section 11.1]{CAV}, we define the Jacobian of \(C\) as \(J(C) = H^0(\omega_C)^*/H_1(C,\Z)\), and for \(\omega = (\omega_1,\dots,\omega_g)\) a basis of \(H^0(\omega_C)\) and the base point \(P_\infty = (0:1:0)\) we define the Abel-Jacobi map
\[\begin{aligned}
\alpha\colon C &\to J(C),\\ 
\quad Q &\mapsto
\int_{P_\infty}^{Q} \omega.
\end{aligned}\]

Choosing a symplectic basis of \(H_1(C,\Z)\) gives rise to the isomorphism
\(J(C) \simeq \C^3/\Omega\Z^3 + \Z^3,\)
where \(\Omega\) is a matrix in the Siegel upper half-space \(\HH_3 = \{ Z\in\C^{3\times 3} : Z = Z^t, \im(Z)>0\}.\) We say that \(\Omega\) is \emph{a (small) period matrix for} \(C\).

The following two classical theorems, due to Riemann and Siegel respectively, deal with the zero locus of the Riemann theta functions and the values of a function of an algebraic curve on non-special divisors. Recall that the \emph{Riemann theta function} \(\theta\colon\C^g \times \HH_g \to \C \) is given by
\[\theta(z, \Omega) = \sum_{n\in\Z^g} \on{exp}(\pi i n^t\Omega n + 2\pi i n^tz),\]
and that a non-special divisor \(D\) is a divisor with \(\ell(K-D) = 0\) for \(K\) a canonical divisor of \(C\).

\begin{theo}[{Riemann's Vanishing Theorem, see \cite[Corollary~3.6]{MumfordI}}]\label{coro:RVT}
	Let~\(C\) be a curve defined over~\(\C\) of genus~\(g\), let \(J(C)\) be the Jacobian of \(C\) with period matrix \(\Omega\in\HH_g\) and let \(\alpha\) be an Abel-Jacobi map of~\(C\).
	There is an element \(\Delta\in J(C)\), called a \emph{Riemann constant} with respect to \(\alpha\), such that the function \(\theta(\,\cdot\,, \Omega)\) vanishes at \(z\in\C^g\) if and only if there exist \(Q_1,\dots,Q_{g-1}\in C\) that satisfy \[\pushQED{\qed} z \equiv \alpha(Q_1 + \dots + Q_{g-1}) - \Delta \mod \Omega\Z^g + \Z^g.\qedhere \popQED\]
\end{theo}
 	
The choice of a base point determines uniquely the Riemann constant \(\Delta\), as shown by Mumford in Theorem~3.10 and Corollary~3.11 of \cite{MumfordI}.

\begin{theo}[{Theorem 11.3 in Siegel \cite{Siegel}}]\label{thm:ImageOfDivisor}
	Let \(C\) be a curve of genus~\(g\) over \(\C\), and let \(\phi\) be a function on \(C\) with
	\[\divv(\phi) = \sum_{i=1}^m A_i - \sum_{i=1}^m B_i.\]
	
	Let \(P\in C\) and let \(\omega\) be a basis of \(H^0(\omega_C)\) for which the Jacobian \(J(C)\) has period matrix~\(\Omega\in\HH_g\). Let \(\Delta\) be the Riemann constant with respect to the Abel-Jacobi map \(\alpha\) with base point \(P\).
	
	Choose paths from the base point \(P\) to \(A_i\) and \(B_i\) that satisfy
	\[\sum_{i=1}^m \int_P^{A_i} \omega = \sum_{i=1}^m \int_P^{B_i} \omega.\]
	
	Then, given an effective non-special divisor \(D = P_1 + \dots + P_g\) of degree~\(g\) that satisfies \(P_j\notin\{A_i,B_i : 1\leq i \leq m\}\),  one has
	\begin{equation}\label{eq:ImageOfDivisor}
	\phi(D):=\phi(P_1) \dots \phi(P_g) = E\prod_{i=1}^m \dfrac{\theta(\sum_{j=1}^g \int_P^{P_j} \omega - \int_P^{A_i} \omega - \Delta, \Omega )}{\theta(\sum_{j=1}^g \int_P^{P_j} \omega - \int_P^{B_i} \omega - \Delta, \Omega )},
	\end{equation}
	where \(E\in\C^\times\) is independent of \(D\), and the integrals from \(P\) to \(P_j\) take the same paths both in the numerator and the denominator. \qed
\end{theo}

Observe that in \eqref{eq:ImageOfDivisor} we are evaluating the Riemann theta functions at points of the Jacobian.

We shall need a version of Theorem~\ref{thm:ImageOfDivisor} in terms of Riemann theta constants. Given \(x = (x_1, x_2)\) with \(x_i\in\R^{g}\), the \emph{Riemann theta constant (with characteristic \(x\))} is the function \(\theta[x]\colon \HH_g \to \C \) given by
\begin{equation}\label{eq:transformation}
\theta \begin{bmatrix}
x_1\\ x_2
\end{bmatrix}(\Omega) = \on{exp}(\pi i x_1^t\Omega x_1 + 2\pi i x_1^tx_2)  \theta(\Omega x_1 + x_2, \Omega) \,.\end{equation}

We use the following two elementary properties of the Riemann theta constants: They are even in \(x\), that is,
\begin{equation}\label{eq:ThetaCharSymmetric}
\theta [x](\Omega) = \theta [-x](\Omega)\,, \end{equation}
and they are quasi-periodic in \(x\), that is, for $m=(m_1,m_2)\in \Z^{2g}$ one has
\begin{equation}\label{eq:ThetaCharQuasiperiodic}
\theta [x+m](\Omega) = \on{exp}(2\pi i x_1 m_2)  \theta [x](\Omega)\,. \end{equation}

Due to the quasi-periodicity of the Riemann theta constants, we must fix representatives in \(\R^{2g}\) for the points of the Jacobian: Throughout, we consider the composition of the maps 
\begin{equation}\label{eq:tilde}
\xymatrix{
	C\ar[r]^-\alpha&J(C) \ar[r]^-{\underline{\cdot}}&\R^{2g}/\Z^{2g}\ar[r]^{\tilde{\cdot}}&[0,1)^{2g}
}
\end{equation}
where \(\alpha\) is the Abel-Jacobi map, the map \(\underline{\cdot}\) identifies \(J(C)\) with \(\R^{2g}/\Z^{2g}\) via \(\Omega x_1 + x_2 \mapsto (x_1, x_2)\) and \(\tilde{\cdot}\) maps a class in \(\R^{2g}/\Z^{2g}\) to its representative with entries in \([0,1)\). For \(P\in C\) we write \(\tilde{P}\) instead of \(\tilde{\underline{\alpha(P)}}\); and in the case of a divisor \(D = \sum n_PP\), we define \(\tilde{D} := \sum n_P\tilde{P} \in\R^{2g}\). Note that with this definition for most divisors \(D\) we get that \(\tilde{D}\) and \(\tilde{\alpha(D)}\) are different.

\medskip With the definitions above, one can rewrite Theorem~\ref{thm:ImageOfDivisor} in terms of Riemann theta constants as follows:

\begin{coro}\label{coro:ImageOfDivisor}
	With the notation of Theorem~\ref{thm:ImageOfDivisor}, let \(a_i\) (resp.~\(b_i\)) be the element in~\(\R^{2g}\) that satisfies \(\int_P^{A_i}\omega = \Omega (a_i)_1 + (a_i)_2\) (resp. \(\int_P^{B_i}\omega = \Omega (b_i)_1 + (b_i)_2\)).
	We have
	\begin{equation*}
	\phi(D) = E'\prod_{i=1}^m \dfrac{\theta\left[\tilde{D} - a_i - \tilde\Delta\right](\Omega)}{\theta\left[\tilde{D} - b_i - \tilde\Delta\right](\Omega )},
	\end{equation*}
	where \(E'\in\C^\times\) is also independent of \(D\).
\end{coro} 

\begin{proof}
	Observe that the exponential factor in \eqref{eq:transformation} for Riemann theta constants can be written as \(\exp(\pi i B(x,x))\) where \(B\) is the symmetric bilinear form given by
	\[ B(u,v) =\, u^t\begin{pmatrix}
	\Omega&\id_g\\\id_g&0
	\end{pmatrix}v.\]
	Let \(Q(u) = B(u,u)\) and let \(c = \tilde{D} - \tilde\Delta\). 
	For \(j = 1, \dots, g\), let \(x_j = \tilde{P_j}\) and choose a path from \(P\) to \(P_j\) that satisfies
	\(\int_P^{P_j}\omega = \Omega (x_j)_1 + (x_j)_2\in \C^g.\)
	
	Let~\(E'\in\C^\times\) be defined by
	\[ E\prod_{i=1}^m \dfrac{\theta\left(\left(\sum_{j=1}^g \int_P^{P_j} \omega\right) - \int_P^{A_i} \omega - \Delta, \Omega \right)}{\theta\left(\left(\sum_{j=1}^g \int_P^{P_j} \omega\right) - \int_P^{B_i} \omega - \Delta, \Omega \right)} = E'\prod_{i=1}^m \dfrac{\theta\left[\tilde{D} - a_i - \tilde\Delta\right](\Omega)}{\theta\left[\tilde{D} - b_i - \tilde\Delta\right](\Omega )}.\]
	
	\bigskip
	We want to prove that \(E'\) does not depend on \(D\). By \eqref{eq:transformation} we get
	\[
	\dfrac{E}{E'} = \exp\left(\pi i \sum_{i=1}^{m}(Q(c - a_i) - Q(c - b_i))\right),
	\]
	so it suffices to show that \(\sum_{i=1}^{m}(Q(c - a_i) - Q(c - b_i))\) does not depend on \(D\). We have
	\begin{equation*}
	\begin{aligned}
	\sum_{i=1}^{m}(Q(c - a_i) - Q(c - b_i)) &= \sum_{i=1}^{m}(Q(a_i) - Q(b_i) - 2B(c, a_i - b_i)) \\
	&= \sum_{i=1}^{m}Q\left(a_i\right) - \sum_{i=1}^{m}Q\left(b_i\right) - 2B\left(c, \sum_{i=1}^{m} (a_i - b_i)\right),
	\end{aligned}
	\end{equation*}
	but we know 
	\[\sum_{i=1}^m \int_P^{A_i} \omega = \sum_{i=1}^m \int_P^{B_i} \omega,\]
	so in terms of characteristics we obtain \(\sum_{i=1}^{m} (a_i - b_i) = 0\) and then it follows that
	\[\sum_{i=1}^{m}(Q(c - a_i) - Q(c - b_i)) = \sum_{i=1}^{m}Q\left(a_i\right) - \sum_{i=1}^{m}Q\left(b_i\right)\]
	does not depend on \(D\).
\end{proof}

\begin{lemma}\label{lem:ImageOfDivisor}
	Let \(C\) be a Picard curve defined over \(\C\) given by $y^3 = x(x-1)(x-\lambda)(x-\mu)$, and consider the branch points 
	$P_0=(0,0)$, $P_1=(1,0)$, $P_\lambda=(\lambda, 0)$, $P_\mu=(\mu,0)$, and $P_\infty$ at infinity. Let \(J(C)\) be the Jacobian of \(C\) with period matrix \(\Omega\), let \(\alpha\) be the Abel-Jacobi map with base point \(P_\infty\), and let \(\Delta\in J(C)\) be the associated Riemann constant.		
	
	Then, for every non-special divisor $D = R_1 + R_2 + R_3$, we have
		\[
	x(D) = E\, \eps(D) \left(\dfrac{\theta[\tilde D - \tilde{P_0} - \tilde\Delta](\Omega )}{\theta[\tilde D - \tilde\Delta](\Omega)}\right)^3,
	\]
	where $\eps(D) = \exp (6\pi i (\tilde D - \tilde{P_0} - \tilde\Delta)_1(\tilde{P_0})_2)$ and $E\in\C^\times$ is a constant independent of $D$.
\end{lemma}

\begin{proof}
	Let \(\omega\) be the basis of holomorphic differentials for which \(J(C)\) has period matrix \(\Omega\). The divisor of the function \(x\) on \(C\) is \(\on{div}(x) = 3\,P_0 - 3\,P_\infty\), so in order to apply Corollary~\ref{coro:ImageOfDivisor} for~\(\phi=x\) and \(P=P_\infty\), we choose three times the zero path from \(P_\infty\) to itself, the path \(\gamma_1\) from \(P_\infty\) to \(P_0\) that for \(a_1 = \tilde{P_0}\) satisfies 
	\[\int_{\gamma_1}\omega = \Omega (a_1)_1 + (a_1)_2\in \C^3,\]
	and paths \(\gamma_2\), \(\gamma_3\) from \(P_\infty\) to \(P_0\) that satisfy
	\begin{equation}\label{eq:Relationms}
	\sum_{k=1}^3\int_{\gamma_k} \omega = 0 \text{ in } \C^3.
	\end{equation}
	Let \(a_2, a_3\) be the elements in \(\R^6\) that satisfy
	\[\int_{\gamma_k} \omega = \Omega (a_k)_1 + (a_k)_2 \text{ for } k = 2,3.\]
	Then, by Corollary~\ref{coro:ImageOfDivisor}, we have
	\begin{equation}\label{eq:prevFromula}
	\begin{aligned}
	\phi(D) = E'\prod_{k=1}^3\dfrac{\theta[\tilde{D} - a_k - \tilde{\Delta}](\Omega )}{\theta[\tilde{D} - \tilde{\Delta}] (\Omega)}
	\end{aligned}
	\end{equation}
	for some constant \(E'\in\C^\times\) independent of \(D\). 
	Note that for \(k=1,2,3\) we have
	\[\underline{P_0} = (a_k \,\on{mod} \Z^{6}),\]
	so the differences \(a_i-a_j\) for \(i\neq j\) are integer vectors. Applying the quasi-periodicity property~\eqref{eq:ThetaCharQuasiperiodic}, equation \eqref{eq:prevFromula} becomes
	\[\phi(D) = E'  \dfrac{\exp(2\pi i (\tilde{D} - \tilde{P_0} - \tilde{\Delta})_1 (a_1 - a_2 + a_1 - a_3)_2)\,\theta[\tilde{D} - \tilde{P_0} - \tilde{\Delta}](\Omega)^3}{\theta[\tilde{D} - \tilde{\Delta}](\Omega)^3}.\]
	But it follows from \eqref{eq:Relationms} that the sum \(a_1 + a_2 + a_3\) is zero, so we obtain \(a_1 - a_2 + a_1 - a_3 = 3a_1 = 3\tilde{P_0}\) and the statement follows.
	\end{proof}
	
	The final step is to choose the right non-special divisors.
	
\begin{theo} \label{thm:thomae-P}
	Let \(C\) be a Picard curve defined over \(\C\) given by $y^3 = x(x-1)(x-\lambda)(x-\mu)$, and consider the branch points 
	$P_0=(0,0)$, $P_1=(1,0)$, $P_\lambda=(\lambda, 0)$, $P_\mu=(\mu,0)$, and $P_\infty$ at infinity. Let \(J(C)\) be the Jacobian of \(C\) with period matrix \(\Omega\), let \(\alpha\) be the Abel-Jacobi map with base point \(P_\infty\), and let \(\Delta\in J(C)\) be the associated Riemann constant.	
	Then, for \(\eta\in\{\lambda, \mu\}\), we have
\begin{equation}\label{eq:Thomae-P}
	\eta = \eps_\eta \left(\dfrac{\theta[\tilde{P_1} + 2\tilde{P_\eta} - \tilde{P_0} - \tilde\Delta](\Omega)}{\theta[2\tilde{P_1} + \tilde{P_\eta} - \tilde{P_0} - \tilde\Delta](\Omega)}\right)^3,
\end{equation}
where \(\eps_\eta = \exp (6\pi i ((\tilde{P_\eta} - \tilde{P_1})_1(\tilde{P_0})_2 + \tilde\Delta_1(3\tilde{P_1} + 3\tilde{P_\eta} - 2\tilde\Delta)_2))\).
\end{theo}

\begin{proof}
We apply Lemma~\ref{lem:ImageOfDivisor} twice, 
to the
divisors $D_1 = P_1 + 2P_\eta$ and $D_2 = 2P_1 + P_\eta$, which Koike--Weng prove that are non-special in \cite[pg.~506]{KW}. Then, we get
\begin{equation}\label{eq:lambdamid}
\begin{aligned}
	\eta = \dfrac{x(P_1)x(P_\eta)^2}{x(P_1)^2x(P_\eta)} & = \dfrac{E'\eps(D_1)\left(\dfrac{\theta[\tilde{P_1} + 2\tilde{P_\eta} - \tilde{P_0} - \tilde\Delta](\Omega)}{\theta[\tilde{P_1} + 2\tilde{P_\eta}) - \tilde\Delta](\Omega)}\right)^3}{E'\eps(D_2)\left(\dfrac{\theta[2\tilde{P_1} + \tilde{P_\eta} - \tilde{P_0} - \tilde\Delta](\Omega)}{\theta[2\tilde{P_1} + \tilde{P_\eta} - \tilde\Delta](\Omega)}\right)^3} \\
	  & = \dfrac{\eps(D_1)}{\eps(D_2)}\left(\dfrac{\theta[\tilde{P_1} + 2\tilde{P_\eta} - \tilde{P_0} - \tilde\Delta](\Omega)}{\theta[\tilde{P_1} + 2\tilde{P_\eta} - \tilde\Delta](\Omega)}\dfrac{\theta[2\tilde{P_1} + \tilde{P_\eta} - \tilde\Delta](\Omega)}{\theta[2\tilde{P_1} + \tilde{P_\eta} - \tilde{P_0} - \tilde\Delta](\Omega)}\right)^3.
\end{aligned}
\end{equation}
Moreover, using the symmetry \eqref{eq:ThetaCharSymmetric} and quasi-periodicity \eqref{eq:ThetaCharQuasiperiodic} of the Riemann theta constants we also obtain
\begin{align*}
\theta[\tilde{D_2} - \tilde\Delta](\Omega) & = \theta[-\tilde{D_2} + \tilde\Delta](\Omega)\\
& = \theta[\tilde{D_1} - \tilde\Delta + \underbrace{2\tilde\Delta - 3\tilde{P_1} - 3\tilde{P_\eta})}_{\in \Z^6}](\Omega)            \\
& = \exp(2\pi i(\tilde{D_1} - \tilde\Delta)_1(2\tilde\Delta - 3\tilde{P_1} -3\tilde{P_\eta})_2))\theta[\tilde{D_1} - \tilde\Delta](\Omega)
\end{align*}
so that \eqref{eq:lambdamid} becomes
\[\eta = \eps_\eta \cdot \left(\dfrac{\theta[\tilde{P_1} + 2\tilde{P_\eta} - \tilde{P_0} - \tilde\Delta](\Omega)}{\theta[2\tilde{P_1} + \tilde{P_\eta} - \tilde{P_0} - \tilde\Delta](\Omega)}\right)^3,\]
with 
\begin{align*}
\eps_\eta & = \dfrac{\eps(D_1)}{\eps(D_2)}\exp(2\pi i(\tilde{D_1} - \tilde\Delta)_1(2\tilde\Delta - 3\tilde{P_1} - 3\tilde{P_\eta})_2)^3 \\
& =  \dfrac{\exp (6\pi i (\tilde{P_1} + 2\tilde{P_\eta} - \tilde{P_0} - \tilde\Delta)_1(\tilde{P_0})_2)}{\exp (6\pi i (2\tilde{P_1} + \tilde{P_\eta}  - \tilde{P_0} - \tilde\Delta)_1(\tilde{P_0})_2)}\exp(6\pi i(\tilde{D_1} - \tilde\Delta)_1(2\tilde\Delta - 3\tilde{P_1} - 3\tilde{P_\eta})_2) \\
& = \exp (6\pi i ((\tilde{P_\eta} - \tilde{P_1})_1(\tilde{P_0})_2 + \tilde\Delta_1(3\tilde{P_1} + 3\tilde{P_\eta} - 2\tilde\Delta)_2 ))
\end{align*}
as desired.
\end{proof}

\begin{rem}\label{rmk:ThomaevsKW}
Compare the above formula in Theorem~\ref{thm:ImageOfDivisor} with the ones given in \cite[Eq. 9]{KW}. The formulas are the same as in \eqref{eq:Thomae-P} replacing \(\eps_\eta\) by 1, hence in general they do not hold due to the absence of the precise root of unity.

However, if we follow the original work by Picard \cite[p. 131]{Pi}, then we obtain a particular form of the period matrix $\Omega$ (see also Shiga \cite[Proposition I-3]{Shg}) for which it is always the case that \(\eps_\lambda = \eps_\mu = 1\). In such case, the formulas in \cite{KW} hold.
\end{rem}

\section{The algorithm\label{sec:algorithm-P}}

In this section we explain how to use the formula in Theorem~\ref{thm:thomae-P} to obtain an inverse Jacobian algorithm for Picard curves, that is, an algorithm that, given the Jacobian of a Picard curve \(C\), returns a model of \(C\). 

The following result characterizes the Jacobian of a Picard curve based on work of Koike-Weng and Estrada.
\begin{prop}\label{prop:charj(P)}
	Let \(X\) be a simple principally polarized abelian variety of dimension 3 defined over an algebraically closed field \(k\). If \(X\) has an automorphism \(\varphi\) of order~3, then we have that \(X\) is the Jacobian of a Picard curve. Furthermore, let \(\rho\) be the curve automorphism \(\rho(x,y) = (x,\zetaCC{3}y)\), and let \(\rho_*\) be the automorphism of the Jacobian that it induces. Then we get \(\langle \varphi \rangle = \langle \rho_* \rangle\).
\end{prop}

\begin{proof}
	By Oort-Ueno \cite{ppavs3}, since \(X\) is a simple principally polarized abelian variety of dimension \(\leq 3\) over an algebraically closed field, then it is the Jacobian of a curve. Let \(C\) be a curve with~\(X\cong J(C)\).
	
	By Torelli's Theorem, see Milne \cite[Section 12]{Milne-JV}, there is some non-trivial automorphism \(\nu\) of \(C\) that satisfies \(\varphi = \pm\nu_*\). Then the automorphism \(\nu^4\), which we call \(\eta\), satisfies \(\eta_* = (\nu^4)_* = (\pm\nu)^4_* = \varphi^4 = \varphi\), hence by the uniqueness in Torelli's Theorem we obtain that \(\eta\) has order~3.
	
	Therefore, the degree of the map~\(\pi\colon C \to C/\langle\eta\rangle\) is also 3, and by the Riemann-Hurwitz formula one obtains that \(C/\langle\eta\rangle\) has either genus 0 or 1.
	But~\(X\) is simple, so the curve \(C/\langle\eta\rangle\) is isomorphic to \(\mathbb{P}^1\) and \(\pi\) has 5 ramification points.
	
	Then \(k(C)/k(C/\langle\eta\rangle)\) is a Kummer extension of degree~3, hence \(C\) is given by an equation of the form \(y^3 = h(x)\) where \(h\) has 4 different roots.
	By Lemma~7.3 in Estrada \cite[Appendix~I]{Hol}, we obtain a model for \(C\) given by \(y^3 = f(x)\) where  \(f\) has degree~4 and distinct roots and \(\eta\) is either the automorphism \(\rho\) given by~\((x,y) \mapsto (x,z_3y)\) or its square.
\end{proof}

\begin{rem}\label{rmk:vsKW}
	While the idea behind the proof is the same in Proposition~\ref{prop:charj(P)} and in \cite[Lemma~1]{KW}, the assumptions in \cite{KW} are in a way more restrictive, as Koike and Weng focus on maximal CM Picard curves. Moreover, the proof in \cite{KW} has a gap, which is fixed exactly by our reference to Estrada \cite[Appendix I]{Hol}.
\end{rem}

It follows from Proposition~\ref{prop:charj(P)} that one can think of the input for this algorithm to be a period matrix \(\Omega\in\HH_3\) together with the rational representation of an automorphism of order 3. To give the curve we will compute the values of \(\lambda\) and \(\mu\) in a Legrendre-Rosenhain equation of the curve.

First we want to determine the points in $\C^3/(\Omega/\Z^3 + \Z^3)$ that correspond to the Riemann constant \(\Delta\) and the image of the branch points via \(\alpha\). The former is given by the following result due to Koike and Weng.

\begin{prop}[{Koike--Weng \cite[Lemma 10]{KW}}] \label{prop:Delta}
	Let \(J(C)\) be the Jacobian of a Picard curve \(C\), let \(\rho_*\) be the automorphism of~\(J(C)\) induced by the curve automorphism \(\rho(x,y) = (x,z_3y)\), and let \(N = \begin{pmatrix}
	\alpha&\beta\\\gamma&\delta
	\end{pmatrix}\in\Sp(6,\Z)\) be the transposed rational representation of \(\rho_*\). Then, the Riemann constant \(\Delta\in J(C)\) is the unique 2-torsion point that satisfies
\[\underline\Delta = (N^{-1})^t{\underline{\Delta}} + \dfrac{1}{2}\begin{pmatrix}
(\gamma^t{\delta})_0\\(\alpha^t{\beta})_0
\end{pmatrix} =:N[\underline\Delta], \]
where \(X_0\) stands for the diagonal of the matrix \(X\).
\end{prop}

The following step is to identify the image under \(\alpha\) of the branch points.

\begin{theo}\label{thm:setDs-P}
	Let \(J(C)\) be the Jacobian of a Picard curve \(C\), let \(\rho_*\) be the automorphism of~\(J(C)\) induced by the curve automorphism \(\rho(x,y) = (x,z_3y)\). Let \(\cB\) be the set of affine branch points of~\(C\), let \(\alpha\) be the Abel-Jacobi map with base point \(P_\infty=(0:1:0)\), let \(\Delta\) be the Riemann constant with respect to \(\alpha\) and define
	\[\Theta_3 := \left\{ x\in J(C)[1-\rho_*] : \theta[x+\underline\Delta](\Omega) = 0\right\}.\]
	
	Then \(\alpha(\cB)\) and \(-\alpha(\cB)\) are the only subsets \(\cT \subset J(C)\) of four elements such that:
	\begin{enumerate}[label=(\roman*)]
		\item[(i)] the sum \(\sum_{x\in \cT} x\) is zero,
		\item[(ii)] \(\cT\) is a set of generators of \(J(C)[1-\rho_*]\), and
		\item[(iii)] the set \(\cO(\cT) :=\{ \sum_{x\in \cT} a_x{x} : a\in\Z^4_{\geq0}, \sum_{x\in \cT} a_x \leq 2 \}\) satisfies \[\cO(\cT) = \Theta_3.\]
	\end{enumerate}
\end{theo}
\begin{proof}
	We first show that \(\alpha(\cB)\) and \(-\alpha(\cB)\) satisfy (i)--(iii), and then we prove that these are the only possibilities.
	
	That \(\alpha(\cB)\) satisfies (i) follows from \(\divv(y) = \sum_{P\in\cB} P - 4P_\infty\). That \(\alpha(\cB)\) satisfies (ii) is proven by Koike and Weng in \cite[Remark 8]{KW}. 
	Next we prove that \(\alpha(\cB)\) satisfies  (iii). On the one hand, given \(Q_1, Q_2\in \cB\cup\{P_\infty\}\) we have \(\alpha({Q_1} + {Q_2}) \in \Theta_3\) by Riemann's Vanishing Theorem~\ref{coro:RVT}, and since we have \(\alpha(P_\infty) = 0\), this implies
	\[\left\{ \sum_{P\in\cB} a_P\alpha(P) : a\in\Z_{\geq 0}^\cB, \sum_{P\in\cB} a_P \leq 2 \right\} \subseteq \Theta_3.\]
	
	To prove the opposite inclusion, let \(x\in\Theta_3\). Since \(x\) satisfies \(\theta[x+\underline{\Delta}](\Omega) = 0\), by Riemann's Vanishing Theorem~\ref{coro:RVT} there exist \(Q_1,Q_2\in C\) such that we have \(x = \alpha(Q_1 + Q_2)\). Moreover, since \(x\) is a~\((1 - \rho_*)\)-torsion point, we get 
	\[\alpha(Q_1 + Q_2) = \rho_*(\alpha(Q_1 + Q_2)) = \alpha(\rho(Q_1)+\rho(Q_2)),\] 
	hence there exists a function \(h\) on \(C\) such that \(\divv(h) = {\rho}(Q_1) + {\rho}(Q_2) - Q_1 - Q_2\). Note now that a Picard curve is non-hyperelliptic, since one checks that the canonical map is the embedding \((x:y:1) \colon C \to \mathbb{P}^2\). Then we conclude that \(h\) is constant, since otherwise it has degree at most~2, hence the curve would be {hyperelliptic}. Therefore we have \({\rho}(Q_1) + {\rho}(Q_2) = Q_1 + Q_2\), but since \(\rho\) has order~3, the cardinality of the orbit of \(Q_i\) has length~3 or~1, we obtain \(\rho(Q_i) = Q_i\). Therefore \(Q_1\) and \(Q_2\) are branch points, so the other inclusion holds.
	
	It is clear that \(-\alpha(\cB)\) satisfies (i) and (ii). To see that it satisfies (iii), it is enough to prove that~\(\Theta_3\) is invariant under the map \(x\mapsto -x\). But this follows from the symmetry of the Riemann theta constants, see \eqref{eq:ThetaCharSymmetric}.
	
	Next we prove that \(\alpha(\cB)\) and \(-\alpha(\cB)\) are, in fact, all the subsets that satisfy~(i)--(iii). 
	
	Let \(B\) denote an ordering of \(\alpha(\cB)\). Given a sequence \(T = (t_1,t_2,t_3,t_4)\) in~\(J(C)^4\) of distinct elements such that the set \( \{t_1, t_2, t_3, t_4\} \) satisfies (i)--(iii), we define the map \(\gamma[T]\colon\FF_3^3\to J(C)[1-{\rho_*}]\) given by \(r \mapsto \sum_{i=1}^{3}r_it_i\). By Remark~8 in Koike--Weng \cite{KW} we have \(\#J(C)[1-\rho_*]\cong(\Z/3\Z)^3\), thus it follows from (i) and (ii) that \(\gamma[T]\) is a bijection. 
	
	Consider the diagram
	\[    \xymatrix{	
		\FF_3^3 \ar[rr]^{M(T)} \ar[rd]_{\gamma[T]} & & \FF_3^3  \ar[ld]^{\gamma[B]}\\
		&J(C)[1-{\rho_*}] &}\]
	where \(M(T)\) is the unique invertible matrix in \(\FF_3^{3\times 3}\) that makes the diagram commutative. Note that choosing a matrix \(M(T)\) determines \(T\) uniquely.
	
	Let \(e_1, e_2, e_3\) be the standard basis vectors of \(\FF_3^3\), and let \(e_4 = -e_1-e_2-e_3\), so for \(i = 1, \dots, 4\) we have \(\gamma[T](e_i) = t_i\). Consider \[\cO_0 = \left\{ \sum_{i=1}^{4} a_ie_i : a\in\Z^4_{\geq 0}, \sum_{i=1}^{4}a_i \leq 2\right\} \subset \FF_3^3.\]
	One can check \(\#\cO_0 = 15\), and moreover we have \(\gamma[T](\cO_0) = \cO(\{t_1, t_2, t_3, t_4\})\). If the set of elements of \(T\) satisfies (iii), then we have
	\begin{equation*}\label{eq:cO maps to Theta_3}
	\gamma[T](\cO_0) = \cO(\{t_1, t_2, t_3, t_4\})= \Theta_3 = \gamma[B](\cO_0),
	\end{equation*}
	and thus \(\cO_0\) is stable under \(M(T)\).
		
	We checked with SageMath \cite{SageMath} that there are exactly 48 invertible matrices in \(\FF_3^{3\times 3}\) that map~\(\cO_0\) to itself. Since a matrix \(M(T)\) determines \(T\) uniquely, there are 48 sequences \(T\in J(C)^4\) that satisfy~(i)--(iii). However, if we vary \(\sigma\) in the symmetric group of 4 letters and \(s\in \{\pm 1\}\), then \(s\sigma(B)\)~gives 48 sequences, which are different.
	We conclude that \(\alpha(\cB)\) and \(-\alpha(\cB)\) are the only subsets of~\(J(C)\) with 4 elements that satisfy (i)--(iii).
\end{proof}

\begin{rem}\label{rmk:findDsvsKW}
With Theorem~\ref{thm:setDs-P}, we make precise the idea hinted in Corollary 11 of Koike--Weng \cite{KW}. There, they claim the existence of a 4-element set that satisfies (i) and (ii), prove that~\(\alpha(\cB)\) does satisfy (i) and (ii), and assume without further comments that when one finds such a set, it is \(\alpha(\cB)\). 

This is problematic not only because they disregard the case where the set is \(-\alpha(\cB)\) but especially because they do not consider (iii) at all, since there exist 4-element sets in \(J(C)\) that satisfy (i) and (ii) which are not \(\alpha(\cB)\) or even \(-\alpha(\cB)\). 

In fact, there are \(\#\on{GL}_3(\FF_3) = 11232\) possible sequences \(T\in J(C)^4\) that satisfy (i) and (ii), hence the probability of finding one that corresponds to a permutation of \(B\) is \(1/468 \approx 0.002\). 
\end{rem}

\bigskip We now have all the tools to state the algorithm.
\begin{alg}{A period matrix \(\Omega\in\HH_3\) of the Jacobian of a Picard curve \(C\), and the transposed rational representation \(N\in\Z^{6\times 6}\) of the automorphism of the Jacobian \(\rho_*\) induced by the curve automorphism \(\rho(x,y) = (x,z_3y)\).}{The complex values \(\lambda\) and \(\mu\) in a Legendre--Rosenhain equation \(y^3 = x(x-1)(x-\lambda)(x-\mu)\) for the Picard curve \(C\).}\label{alg:Picard}
	\item  Let \(D\) be the unique solution of \(N[D] = D\) in \(\frac{1}{2}\Z^6/\Z^6\).
	\item  Compute the set
	\[\underline{\Theta_3} = \left\{x\in\frac{1}{3}\Z^6/\Z^6 : N^tx = x \text{ and } \theta[x + D](\Omega) = 0\right\}\]
	of cardinality 15.
	\item  Let \(T = \{t_1, t_2, t_3, t_4\}\subset \underline{\Theta_3}\) be a 4-element set that satisfies 
	\begin{enumerate}[label=\roman*.]
		\item \(\sum_{i=1}^4 t = 0\),
		\item \( \{t_1, t_2, t_3\}\) are linearly independent over \(\Z/3\Z\), and
		\item \(\{ \sum_{i=1}^4 a_i{t_i} : (a_i)_i\in\Z^4_{\geq0}, \sum_{i=1}^4 a_i \leq 3 \} = \underline{\Theta_3}\).
	\end{enumerate}
	\item  Compute
	\[\eps_\lambda =  \exp (6\pi i ((\tilde{t_3} - \tilde{t_2})_1(\tilde{t_1})_2 + (\tilde{t_2} + 2\tilde{t_3} - \tilde{D})_1(2\tilde{D} - 3(\tilde{t_2} + \tilde{t_3}))_2 )),\]
	\[\eps_\mu =  \exp (6\pi i ((\tilde{t_4} - \tilde{t_2})_1(\tilde{t_1})_2 + (\tilde{t_2} + 2\tilde{t_4} - \tilde{D})_1(2\tilde{D} - 3(\tilde{t_2} + \tilde{t_4}))_2 )),\]
	and
	\[\lambda = \eps_\lambda \left(\dfrac{\theta[\tilde{t_2} + 2\tilde{t_3} - \tilde{t_1} - \tilde{D}](\Omega)}{\theta[2\tilde{t_2} + \tilde{t_3} - \tilde{t_1} - \tilde{D}](\Omega)}\right)^3,\]
	\[\mu = \eps_\mu \left(\dfrac{\theta[\tilde{t_2} + 2\tilde{t_4} - \tilde{t_1} - \tilde{D}](\Omega)}{\theta[2\tilde{t_2} + \tilde{t_4} - \tilde{t_1} - \tilde{D}](\Omega)}\right)^3.\]
	\item  Return \(\lambda\) and \(\mu\).
\end{alg}

\begin{rem}
	Algorithm~\ref{alg:Picard} is a \emph{mathematical} algorithm, but, because it involves infinite sums, complex numbers and exponentials, it cannot be run on a Turing machine or a physical computer. To do so one needs to truncate the sum on the Riemann theta constants, approximate complex numbers and keep track of the error propagation. For implementation details, we refer the reader to \cite[Section~1.5]{thesis}.
\end{rem}

\begin{proof}[Proof of Algorithm~\ref{alg:Picard}]
	Let \(\Delta\in J(C)\) be the Riemann constant with respect to \(P_\infty = (0:1:0)\) and let \(\cB\) be the set of affine branch points of \(C\).
	By Proposition~\ref{prop:Delta}, the point \(\Delta\) is the only one that satisfies \(N[\underline\Delta] = \underline\Delta\) and is a 2-torsion point, that is, it satisfies \(\underline{\Delta}\in\frac{1}{2}\Z^6/\Z^6\). We conclude~\(D = \underline{\Delta}\).
	
	By Theorem~\ref{thm:setDs-P}, the sequence \((t_1,t_2,t_3,t_4)\) is an ordering of either \(\alpha(\cB)\) or \(-\alpha(\cB)\). In the former case, the values \(\lambda\), \(\mu\) obtained in Step 4 are the \(x\)-coordinates of the affine branch points different from \((0,0)\) and \((0,1)\). A quasi-periodicity argument similar to those in the proofs of Lemma~\ref{lem:ImageOfDivisor} or Theorem~\ref{thm:thomae-P} yields that in the latter case the same holds too.
\end{proof}

\section{Implementation details and CM examples}\label{Llista}

Assume that a Picard curve \(C\) has a model \(y^3 = h(x)\), with \(h(x)\) a polynomial over a number field. After numerically approximating the \(x\)-coordinates of the branch points of \(C\) with Algorithm~\ref{alg:Picard}, we obtain a polynomial \[f(x) = x(x-1)(x-\lambda)(x-\mu)\in\C[x]\] up to some precision, which gives an approximate model for the curve we seek.

Given the quartic polynomial \[p(x) = x^4 + g_2x^2+g_3x+g_4\text{ with } g_2\neq 0\]
we define the \emph{absolute invariants of }\(p(x)\) as 
\[j_1 = \dfrac{g_3^2}{g_2^3}, \qquad j_2 = \dfrac{g_4}{g_2^2}.\]

In order to find \(h(x)\) from \(f(x)\) (when possible), we compute the {absolute invariants} of \(C\) by computing \(j_1\) and \(j_2\) for our approximation of the curve \(C\). We then recognize \(j_1\) and \(j_2\) as algebraic numbers and reconstruct \(h(x)\) from the exact absolute invariants, obtaining 
\[y^3 = h(x) = x^4 + j_1x^2 + j_1^2x + j_1^2j_2.\] 

Note that in order to  be able to recognize \(j_1\) and \(j_2\) as algebraic numbers we have to compute~\(\lambda\) and \(\mu\) with enough precision. 

\bigskip 
One possible application for Algorithm~\ref{alg:Picard} is to compute \emph{maximal CM Picard curves}, that is, Picard curves such that their Jacobians have an endomorphism ring isomorphic to the maximal order of a sextic CM-field~\(K\). Since \(\rho_*\) is an automorphism of order 3, the field \(K\) contains a primitive 3rd~root of unity \(\zeta_3\in K\). {In fact, the field \(K\) is determined by a totally real cubic field~\(K_0\) that satisfies \(K =  K_0(\zeta_3)\).}

Van Wamelen \cite{VanW} gives an algorithm that, given a CM-field \(K\), lists all the isomorphism classes of period matrices of principally polarized abelian varieties with complex multiplication by \(\cO_K\). This method is based on the CM theory due to Shimura and Taniyama, see \cite{shimura-taniyama}. 

If we apply said method to a sextic CM-field containing a primitive third root of unity \(\zeta_3\in K\), then we obtain a list of period matrices corresponding to principally polarized abelian threefolds with CM by \(\cO_K\) with an order-3 automorphism associated to \(\zeta_3\) which, by Proposition~\ref{prop:charj(P)}, are Jacobians of Picard curves. To then obtain the rational representation of the order-3 automorphism is a matter of keeping track of the changes of basis throughout van Wamelen's method, which completes the input for our algorithm. 

\medskip Using Algorithm~\ref{alg:Picard} on the resulting list of pairs \((\Omega,N)\), we computed heuristic models of some maximal CM Picard curves. In particular, the list below contains all maximal CM Picard curves whose CM-field has class number $h\leq 4$. We get the sextic fields from \cite[Table 3]{Klist}, where the complete list of all imaginary abelian sextic number fields with class number $h\leq 11$ is given. 

It follows from K{\i}l{\i}cer \cite[Theorem 4.3.1]{PinarThesis} that our list also includes conjectural models for all Picard curves defined over \(\Q\) with maximal CM over \(\C\), see also \cite[Table 3.1]{PinarThesis}. The curves (1)--(5) also appear in \cite[Section 6.1]{KW}. 

\begin{enumerate}[itemsep=3pt, label=(\arabic*), wide, labelwidth=!, labelindent=0pt]
	\item $y^3 = x^4 - x$, with $K_0$ defined by 
	$\nu^3 - 3\nu -1$.
	\item $y^3 = x^4 - 2\cdot7^2 \, x^2 + 2^3\cdot7^2\, x - 7^3$, with $K_0$ defined by $\nu^3 - \nu^2 - 2\nu +1$.
	\item $y^3 = x^4 - 2\cdot 7^2\cdot 13\, x^2 +2^3\cdot 5\cdot 13\cdot 47\, x - 5^2\cdot 13^2\cdot 31$, with $K_0$ defined by $\nu^3 - \nu^2 - 4\nu  -1$.
	\item $y^3 = x^4 - 2\cdot 7\cdot 31\cdot 73\, x^2 + 2^{11}\cdot 31\cdot 47\, x - 7\cdot 31^2\cdot 11593$, with $K_0$ defined by $\nu^3 + \nu^2 - 10\nu  -8$.
	\item $y^3 = x^4 -2\cdot 7\cdot 43^2\cdot 223\, x^2 + 2^7\cdot 11\cdot 41\cdot 43^2\cdot 59\, x- 11^2\cdot 43^3\cdot 419\cdot 431$, with $K_0$ defined by~$\nu^3 - \nu^2 - 14\nu  -8$.
	\item $y^3 = x^4 - 2 \cdot 3^{2} \cdot 5^{2} \cdot 7^{2}\, x^2 + 2^{9} \cdot 7^{2} \cdot 71\, x - 3^{2} \cdot 5 \cdot 7^{3} \cdot 2621 $, with 
	$K_0$ defined by $\nu^3 - 21\nu  -28$.
	\item $y^3 = x^4 - 2^{2} \cdot 3^{2} \cdot 7^{2} \cdot 37\, x^2 + 5 \cdot 7^{2} \cdot 149 \cdot 257 \, x - 2 \cdot 3^{2} \cdot 5^{2} \cdot 7^{3} \cdot 2683$, with $K_0$ defined by $\nu^3 - 21\nu  +35$.
	\item $y^3 = x^4 -2 \cdot 3^{2} \cdot 5^{2} \cdot 7 \cdot 11 \cdot 13\, x^2 + 2^{7} \cdot 11 \cdot 13 \cdot 59 \cdot 149 \, x - 3^{2} \cdot 5 \cdot 7 \cdot 13^{2} \cdot 17 \cdot 17669$, with $K_0$ defined by $\nu^3 - 39\nu  +26$.
	\item For $K_0$ defined by $\nu^3 - \nu^2 - 6\nu  + 7$,  and $w^3 = 19$,
	$$y^3 = x^4 + (10w^2 - 2w-70)\, x^2 + (96w^2 - 7w-496)\, x + (235w^2 - 215w - 1101).$$
	
	\item For $K_0$ defined by $\nu^3 - \nu^2 - 12\nu  - 11$,  and $w^3 = 37$,
	\begin{equation*}
	\begin{aligned}
	y^3 = \, \, & x^4 + (-2366w^2 + 490w + 24626) \, x^2 + (-257958w^2 - 686928w\\ & + 5152928)\, x
	+ (1226851w^2 - 56922233w + 176054907).
	\end{aligned}
	\end{equation*}
	
	\item For $K_0$ defined by $\nu^3 - 109\nu  - 436$, and $w^3 = 109$,
	\begin{equation*}
	\begin{aligned}
	y^{3} =\, \,  & x^{4}  + \left(1115888872 w^{2} - 4007074778 w - 6321528472\right) x^{2}     \\& + \left(-39141169182336 w^{2} + 294349080537984 w - 512926132238464\right) x \\
	& + 816342009554519305 w^{2} - 9276324622428605048 w\\& + 25684086855493144296.
	\end{aligned}
	\end{equation*}
	
	\item For $K_0$ defined by $\nu^3 - \nu^2 - 42\nu  - 80$, and $w^3 = 127$,	
	\begin{equation*}
	\begin{aligned}
	y^{3} =&\;  x^{4}  + \left(-92075757704 w^{2} + 319193013538 w + 721950578888\right) x^{2}\\ &+ \big(-49404281036538240 w^{2} - 182817463505393280 w + 
	\\& 2167183294305193600\big) x + 21690511027003736433025 w^{2} - 
	\\&118803029086722205449800 w + 49134882128483485627800.
	\end{aligned}
	\end{equation*}
	\item For $K_0$ defined by \(v^3 - 61v - 183\), we have four curves. The first one is defined over $\Q$.
	\[
	\begin{aligned}
	y^3 =\, \, & x^{4} - 2 \cdot 3 \cdot 7 \cdot 61^{2} \cdot 1289 \, x^{2} + 2^{3} \cdot 3^{7} \cdot 11 \cdot 41 \cdot 53 \cdot 61^{2} \, x\\ & - 3^{2} \cdot 7 \cdot 11^{2} \cdot 61^{3} \cdot 419 \cdot 4663\\[8pt]
	y^{3} = \, \, & x^{4} + \left(89264 v^{2} - 547484 v - 4059720\right) x^{2} + \big(-29558196 v^{2} + 49526073 v \\& + 772138494\big) x + 88325678 v^{2} - 16281030326 v - 72348132021
	\end{aligned}
	\]

	\item For $K_0$ defined by $v^3 - v^2 - 22v - 5$,  similarly one gets:
	\[
	\begin{aligned}
	y^3 = \, \, & x^{4}  + 2 \cdot 7 \cdot 67 \cdot 179 \, x^{2} + 2^3 \cdot 3^3 \cdot 5 \cdot 67 \cdot 137 \, x + 5^2 \cdot 7 \cdot 67^2 \cdot 71 \cdot 89\\[8pt]
	y^{3} = \, \, & x^{4}  + \left(12222 v^{2} - 263088 v - 1290744\right) x^{2} + \big(-19721880 v^{2} + 232016400 v\\& + 1277237160\big) x + 11453819175 v^{2} - 62791404525 v - 447679991475 \,.
	\end{aligned}
	\]
	
\end{enumerate}
\begin{rem}	We observe that for each one of the curves listed above, the corresponding unramified field $K(j_1, j_2)$ over $K$ coincides
with the Hilbert class field of $K$, except for the cases $h=h^*=4$ where $K(j_1, j_2) = K$ (see \cite[Main Theorem 1]{ShiTa}).
\end{rem}

\bibliographystyle{unsrt}

\appendix
\newpage
\section{Appendix (by Christelle Vincent)}

Let $C$ be a hyperelliptic curve of genus $g$ defined over $\mathbb{C}$, and denote by $x \colon C \to \mathbb{P}^1$ a morphism of degree 2 from $C$ to $\mathbb{P}^1$. Then $x$ has $2g+2$ branch points which do not depend on the choice of $x$. We fix once and for all an ordering of these branch points, and denote them by $P_1, P_2, \ldots, P_{2g+2}$. Furthermore, for simplicity of notation in what follows we will denote
\begin{equation}
a_j = x(P_j).
\end{equation}

The significance of these quantities is the following: If $x(P_j) \neq \infty$ for any $j$, then a model for $C$ over~$\mathbb{C}$ is given by
\begin{equation}
y^2 = \prod_{j=1}^{2g+2} (x-a_j),
\end{equation}
whereas if there is $k$ with $x(P_k) = \infty$, a model for $C$ over $\mathbb{C}$ is given by
\begin{equation}
y^2 = \prod_{j \neq k} (x-a_j).
\end{equation}

Our goal in this Appendix is to show the following Proposition, which generalizes a formula given by Takase \cite[Theorem 1.1]{takase}. In the statement we use the notation $[a_l,a_m,a_k,a_{\infty}]$ for the cross-ratio
\begin{equation}
[a_l,a_m,a_k,a_{\infty}] = \frac{a_k-a_l}{a_k-a_m} \cdot \frac{a_{\infty}-a_m}{a_{\infty}-a_l}.
\end{equation}
\begin{prop}\label{prop:MainProp}
	Let $C$ be a hyperelliptic curve defined over $\mathbb{C}$, $x \colon C \to \mathbb{P}^1$ be a morphism of degree~$2$ with branch points $P_1, \ldots, P_{2g+2}$, and $\Omega$ be a (small) period matrix for $J(C)$, the Jacobian of $C$. Let $k$, $l$ and $m$ be distinct and belong to the set $\{1,2, \ldots, 2g+2\}$, and fix $P_\infty$ a distinguished branch point of $x$, $\infty \neq k,l,m$. Then, for $a_j = x(P_j)$ and $\eta$ an eta-map associated to $\Omega$ and the base point $P_\infty$ (see Section \ref{sec:preliminaries} for more on eta-maps) with corresponding $U$-set $U_\eta$, we have
	\begin{equation}
	[a_l,a_m,a_k,a_{\infty}]= \exp(4 \pi i (\eta_{m}-\eta_{l})_1(\eta_k)_2)\left( \frac{\theta[\eta_{U_\eta\circ(V \cup \{k,l\})}](\Omega)\theta[\eta_{U_{\eta} \circ (W \cup \{k,l\})}](\Omega)}{\theta[\eta_{U_\eta \circ (V \cup \{k,m\})}](\Omega)\theta[\eta_{U_{\eta} \circ (W \cup \{k,m\})}](\Omega)}\right)^2,
	\end{equation}
	where $V$ and $W$ are any sets that give a disjoint decomposition 
	\begin{equation}
	\{1,2,\ldots, 2g+1,2g+2\} = V \cup W \cup \{k,l,m,\infty\},
	\end{equation}
	with $\#V = \#W = g-1$.
\end{prop}

Indeed, in his work Takase gives the formula above, but only for certain choices of period matrix $\Omega$ for the Jacobian of $C$, which are those given by Mumford \cite{mumford}, using his ``traditional" choice of symplectic basis for the first homology group of the Jacobian. Following this, our earlier article \cite[Theorem 4.5]{BILV} claimed to give the formula for all period matrices, but there remained a mistake in the sign, which had not been corrected to account for the general case. The formula we finally give here is valid for all period matrices, and gives the correct sign. We note that the software available at \cite{BILVcode} has been updated to be correct. We also note that our formula does not assume that $a_{\infty} = \infty$, which explains why we compute the cross-ratio $[a_l,a_m,a_k,a_{\infty}]$ rather than the simpler quotient $\frac{a_k-a_l}{a_k-a_m}$.

As an immediate Corollary, if we denote by $\lambda_i$ for $i = 3, 4, \ldots, 2g+1$ the \emph{Rosenhain invariants} of~$C$, by which we mean the constants appearing in a choice of Rosenhain model
\begin{equation}
C : y^2 = x(x-1) \prod_{i = 3}^{2g+1} (x- \lambda_i)
\end{equation}
for the curve $C$, we obtain the following formula:

\begin{coro}\label{cor:Rosenhain}
	Let $C$ be a hyperelliptic curve defined over $\mathbb{C}$, and fix a choice of Rosenhain model for $C$. Let $P_\infty$ denote the point of $C$ that is ``at infinity" in the Rosenhain model of $C$, $\Omega$ be a choice of period matrix for $J(C)$, the Jacobian of $C$, and $\eta$ be an eta-map associated to $\Omega$ and the base point $P_\infty$ with corresponding $U$-set $U_\eta$. Then for $j \in \{3,4, \ldots, 2g+1\}$, the Rosenhain invariants of $C$ are given by the expression
	\begin{equation}
	\lambda_j =  \exp(4 \pi i (\eta_j-\eta_2)_1(\eta_1)_2)\left( \frac{\theta[\eta_{U_\eta\circ(V \cup \{1,2\})}](\Omega)\theta[\eta_{U_{\eta} \circ (W \cup \{1,2\})}](\Omega)}{\theta[\eta_{U_\eta \circ (V \cup \{1,j\})}](\Omega)\theta[\eta_{U_{\eta} \circ (W \cup \{1,j\})}](\Omega)}\right)^2,
	\end{equation}
	where $V$ and $W$ are two sets of cardinality $g-1$ such that
	\begin{equation}
	V \cup W = \{3, 4, \ldots, 2g+1\} \setminus \{j\},
	\end{equation}
	and the notation $\circ$ denotes the symmetric difference of two sets: For $S,T \subseteq \{1,2, \ldots, 2g+2 \}$, we have
	\begin{equation}\label{eq:circ}
	S\circ T = (S\cup T) \setminus (S \cap T).
	\end{equation}
\end{coro}

\subsection*{Acknowledgements}
The author of this Appendix wishes to thank first and foremost Sorina Ionica, who verified the result independently with a proof that follows Takase's work more closely. She also thanks Marco Streng for first bringing to her attention the need to generalize Takase's work, and Anna Somoza for pointing out that her work on Picard curves could be adapted to obtain the correct sign. Finally, she thanks Jeroen Sijsling for performing computations confirming the correctness of the sign as computed in this Appendix.

\subsection{Preliminaries}\label{sec:preliminaries}

Following the technique used in the body of the paper, we will use Siegel's Theorem \ref{thm:ImageOfDivisor} applied to a suitable choice of function $\phi \colon C \to \mathbb{P}^1$ to obtain our results. To apply Siegel's Theorem, we first need a non-special divisor on $C$:

\begin{lemma}\label{lem:Nonspecial}
	Let $C$ be a hyperelliptic curve defined over $\mathbb{C}$, $x \colon C \to \mathbb{P}^1$ be a morphism of degree 2 from $C$ to $\mathbb{P}^1$, and $P_1,\ldots P_{2g+2}$ be the branch points of $x$. Let $I \subset \{1, 2, \ldots, 2g+2\}$ be any subset of cardinality $g$. Then 
	\begin{equation}
	D = \sum_{i \in I} P_i
	\end{equation}
	is a non-special divisor on $C$. In other words, any sum of $g$ distinct branch points of $x$ is a non-special divisor on $C$.
\end{lemma}

\begin{proof}
	We recall that a divisor $D$ is non-special if $\ell(K-D) =0$, where $K$ is a canonical divisor on the curve $C$. By Riemann-Roch we have that
	\begin{equation}
	\ell(D) - \ell(K-D) = \deg(D) - g+1,
	\end{equation}
	and here $\deg(D) = g$, so to show that $\ell(K-D) =0$ it suffices to show that $\ell(D) = 1$.
	
	Let $D$ be as in the statement of the Lemma, and let $P_{\infty}$ be a branch point of $x$ that does not belong to the support of $D$. For $i \in I$, the function
	\begin{equation}
	x_i(P) = \frac{x(P)-x(P_i)}{x(P)-x(P_\infty)}
	\end{equation}
	has a double zero at $P_i$ and a double pole at $P_\infty$. As a result, the divisor $D$ is equivalent to the divisor
	\begin{equation}
	2gP_\infty - D,
	\end{equation}
	and the linear space associated to the divisor $D$ is isomorphic to the linear space associated to the divisor $2gP_\infty - D$. In particular, their dimensions are the same. Therefore we may show that
	\begin{equation}
	\ell(2gP_\infty - D) = 1
	\end{equation}
	to prove our claim.
	
	Now the linear space associated to the divisor $2gP_\infty - D$ is the space of functions with a pole of order at most $2g$ at $P_\infty$ and zeroes at each of the $g$ points that belong to the support of $D$. We can give a model
	\begin{equation}
	s^2 = f(t)
	\end{equation}
	for our curve $C$, where $f$ is a polynomial of degree $2g+1$ with zeroes at $x(P_j)$, $P_j$ a branch point of $x$, $P_j \neq P_\infty$. Then the function field of $C$ is generated over $\mathbb{C}$ by $s$ and $t$.  Furthermore, we have that
	\begin{equation}
	\divv(s) = \sum_{P_j \neq P_\infty} P_j - (2g+1) P_\infty
	\end{equation}
	and $t$ has a double pole at $P_\infty$ and no other poles. From this it follows that the space of functions with a pole of order at most $2g$ at $P_\infty$ and no other poles is the space of polynomials in $t$ of degree at most $g$. If we require further that this function vanishes at the $g$ points in the support of $D$, we obtain a space of dimension $1$, which completes the proof.
\end{proof}

Secondly, to connect our result to the established literature on hyperelliptic curves, we will need an eta-map associated to a period matrix $\Omega$ and a base point $P_\infty$. We refer the interested reader to either Poor's work \cite{poor} or our earlier work \cite{BILV} for more details on these maps, and present here only the barest of facts necessary to keep this Appendix readable. Let $P_{\infty}$ be an arbitrary but fixed branched point of the degree $2$ morphism $x \colon C \to \mathbb{P}^1$ fixed above, and recall that we have labeled the branch points of $x$ to be $P_1, P_2, \ldots, P_{2g+2}$ (one of these is of course also labeled $P_\infty$). As in the body of the paper, fix $\alpha$ an Abel-Jacobi map for $C$ with base point $P_{\infty}$. Then for $j \in \{1,2, \ldots, 2g+2\}$, we write
\begin{equation}\label{eq:halfint}
\eta_j = \widetilde{P_j} \in \left\{0,\frac{1}{2}\right\}^{2g}
\end{equation}
where $\tilde{\cdot}$ is the map given in equation \eqref{eq:tilde}, and as in the body of the paper we denote the composition of the three maps by the last. (The fact that the coordinates of $\eta_j$ for each $j$ are half-integers follows from the fact that $P_j-P_{\infty}$ is two-torsion in $J(C)$, see \cite[Corollary 2.11]{mumford}.) Furthermore, for any subset $S \subseteq \{1,2, \ldots, 2g+2\}$, we write
\begin{equation}
\eta_S = \sum_{j \in S} \eta_j.
\end{equation}
Note that we use the same convention as in the body of the paper regarding summation. It then follows that
\begin{equation}
\eta_S = \widetilde{D_S},
\end{equation}
for
\begin{equation}
D_S = \sum_{j \in S} P_j.
\end{equation}
We note that the dependence of the eta-map on the period matrix $\Omega$ happens explicitly via the map $\underline{\cdot}$.

Under these assumptions, there exists a subset $U_{\eta}\subseteq\{1,2,\ldots,2g+2\}$ such that
\begin{equation}\label{eq:Uset}
\eta_{U_{\eta}} \equiv \widetilde{\Delta} \pmod{\mathbb{Z}^{2g}}
\end{equation} 
where $\Delta$ is the Riemann constant associated to the choice of Abel-Jacobi map $\alpha$ that we made. We note that in fact there are several such sets; it is customary to choose one of even cardinality, and we have adopted in earlier work the convention that $U_{\eta}$ should also contain $\infty$. This determines the set $U_{\eta}$ uniquely. We call this set a \emph{$U$-set corresponding to $\eta$}. Finally, one can show that if $S$ is the complement of $T$ inside of $\{1,2,\ldots, 2g+2\}$, then 
\begin{equation}\label{eq:complements}
\eta_S = \eta_T.
\end{equation}

\subsection{Proof of the formula}

With this notation and preliminaries in place, we may begin the proof. We begin with an auxiliary result:

\begin{lemma}\label{lem:twopaths}
	Let $P_j$ and $P_{\infty}$ be two distinct branch points of the morphism $x$, $\alpha$ be an Abel-Jacobi map with base point $P_\infty$, and $\gamma$ be a path from $P_{\infty}$ to $P_j$ such that if $\widetilde{P}_j = \eta_j$ (where the map $\tilde{\cdot}$ is as in equation \eqref{eq:tilde}), then
	\begin{equation}
	\int_{\gamma} \omega = \Omega (\eta_j)_1 + (\eta_j)_2.
	\end{equation}
	In this case there exists a second path $\tilde{\gamma}$ from $P_{\infty}$ to $P_j$ such that
	\begin{equation}
	\int_{\gamma} \omega + \int_{\tilde{\gamma}} \omega = 0 \quad \text{in} \quad \mathbb{C}^g.
	\end{equation}
\end{lemma}

\begin{proof}
	We have that $\widetilde{P}_j = \eta_j \in \{0,\frac{1}{2}\}^{2g}$ (see equation \eqref{eq:halfint} and the discussion surrounding it for this fact). From this it follows that if $L_{\Omega} = \Omega \mathbb{Z}^{2g} +\mathbb{Z}^{2g}$ is the lattice attached to the period matrix $\Omega$, we have that
	\begin{equation}
	\int_{\gamma} \omega \in \frac{1}{2}L_{\Omega},
	\end{equation}
	or
	\begin{equation}
	2\int_{\gamma} \omega \in L_{\Omega}.
	\end{equation}
	As a consequence,  $\int_{\gamma} \omega$ and $-\int_{\gamma} \omega$ differ by an element of $L_{\Omega}$, and since every $L_{\Omega}$-translate of $\int_{\gamma} \omega$ is $\int_{\tilde{\gamma}} \omega$ for some other path $\tilde{\gamma}$ from $P_{\infty}$ to $P_j$, it follows that there is $\tilde{\gamma}$ from $P_{\infty}$ to $P_j$ such that
	\begin{equation}
	-\int_{\gamma} \omega = \int_{\tilde{\gamma}} \omega.
	\end{equation}
\end{proof}

We can now give the crucial part of the proof:

\begin{lemma}\label{lem:mainlemma}
	Let $C$ be a hyperelliptic curve defined over $\mathbb{C}$, $x \colon C \to \mathbb{P}^1$ be a morphism of degree $2$ with branch points $P_1, \ldots, P_{2g+2}$, and $\Omega$ be a period matrix for $J(C)$, the Jacobian of $C$. Let $k$, $l$ and $m$ be distinct and belong to the set $\{1,2, \ldots, 2g+2\}$, and fix $P_\infty$ a distinguished branch point of $x$, with $\infty \neq k,l,m$. Then, for $a_j = x(P_j)$, and $\eta$ an eta-map associated to $\Omega$ and to the base point $P_\infty$ with corresponding $U$-set $U_\eta$, we have
	\begin{equation}
	[a_l,a_m,a_k,a_{\infty}]=\epsilon(k,l,m)\left( \frac{\theta[\eta_{S_l \circ U_\eta}](\Omega)\theta[\eta_{T_m \circ U_\eta}](\Omega)}{\theta[\eta_{S_m \circ U_\eta}](\Omega)\theta[\eta_{T_l \circ U_\eta}](\Omega)}\right)^2,
	\end{equation}
	where
	\begin{equation}
	\epsilon(k,l,m) =  \exp(4 \pi i (\eta_{m}-\eta_{l})_1(\eta_k)_2),
	\end{equation}
	and for $j = l,m$, we have
	\begin{equation}
	T_j = V \cup \{j\},
	\end{equation}
	and
	\begin{equation}
	S_j = T_j \cup \{k\} = V\cup\{j,k\},
	\end{equation}
	where $V$ is any set of cardinality $g-1$ such that $V \subset \{1,2,\ldots,2g+2\}$, $k,l, m, \infty \not \in V$.
\end{lemma}

\begin{proof}
	To begin, fix $\infty \in \{1,2,\ldots,2g+2\}$, $\infty \neq k,l,m$, and let
	\begin{equation}
	x_k(P) \colon C \to \mathbb{P}^1
	\end{equation}
	be given by
	\begin{equation}
	x_k(P) = \frac{x(P) - x(P_k)}{x(P)-x(P_{\infty})}.
	\end{equation}
	Then the cross-ratio we seek is given by 
	\begin{equation}
	[a_l,a_m,a_k,a_{\infty}]=  \frac{x_k(P_l)}{x_k(P_m)}.
	\end{equation}
	Next we fix a subset $V \subset \{1,2,\ldots,2g+2\}$ of cardinality $g-1$ such that $k,l, m, \infty \not \in V$. (Note that this is possible since $2g-2\geq g-1$ for $g \geq 1$.) Then the quantity which interests us is given by
	\begin{equation}
	[a_l,a_m,a_k,a_{\infty}]= \frac{x_k(P_l) \prod_{i \in V} x_k(P_i)}{x_k(P_m)\prod_{i \in V} x_k(P_i)}.
	\end{equation}
	In addition, for $j = l,m$, the divisor
	\begin{equation}
	D_j = P_j + \sum_{i \in V} P_i
	\end{equation}
	is a sum of $g$ distinct branch points of $x$, and therefore an effective non-special divisor by Lemma \ref{lem:Nonspecial}.

	Using the notation of Siegel's Theorem~\ref{thm:ImageOfDivisor}, we have
	\begin{equation}
	[a_l,a_m,a_k,a_{\infty}]= \frac{x_k(D_l)}{x_k(D_m)},
	\end{equation}
	and now wish to apply Corollary~\ref{coro:ImageOfDivisor} to compute the quantities $x_k(D_l)$ and $x_k(D_m)$. To do so, we note that
	\begin{equation}
	\divv(x_k) = 2P_k - 2P_{\infty}
	\end{equation}
	and that the supports of the divisors $D_l$ and $D_m$ avoid the support of $\divv(x_k)$. As in the previous section, we denote by $\Delta$ the Riemann constant for the Abel-Jacobi map $\alpha$ of $C$ with base point $P_{\infty}$. In the application of Siegel's Theorem, we will choose the paths from $P_{\infty}$ to $P_{\infty}$ to be the trivial paths. As in Lemma \ref{lem:twopaths}, we fix a path $\gamma$ from $P_k$ to $P_{\infty}$ such that
	\begin{equation}
	\int_{\gamma} \omega = \widetilde{P_k} = \Omega(\eta_k)_1+(\eta_k)_2,
	\end{equation}
	and denote by $\tilde{\gamma}$ the path from $P_k$ to $P_{\infty}$ such that
	\begin{equation}
	\int_{\gamma} \omega + \int_{\tilde{\gamma}}\omega = 0.
	\end{equation}
	We have then that
	\begin{equation}
	\int_{\tilde{\gamma}}\omega = -\widetilde{P}_k.
	\end{equation}
	
	Then if we apply Corollary~\ref{coro:ImageOfDivisor} to $x_k(D_l)$ and $x_k(D_m)$, we obtain
	\begin{align}
	[a_l,a_m,a_k,a_{\infty}] &= \frac{x_k(D_l)}{x_k(D_m)} \\
	& = \left(E' \frac{\theta[\widetilde{P}_l + \sum_{i \in V} \widetilde{P}_i-\widetilde{P}_k-\tilde{\Delta}](\Omega)\theta[\widetilde{P}_l + \sum_{i \in V} \widetilde{P}_i+\widetilde{P}_k-\tilde{\Delta}](\Omega)}{\theta[\widetilde{P}_l + \sum_{i \in V} \widetilde{P}_i-\tilde{\Delta}](\Omega)^2}\right) \\  \notag
	& \quad \div \left(E' \frac{\theta[\widetilde{P}_m + \sum_{i \in V} \widetilde{P}_i-\widetilde{P}_k-\tilde{\Delta}](\Omega)\theta[\widetilde{P}_m + \sum_{i \in V} \widetilde{P}_i+\widetilde{P}_k-\tilde{\Delta}](\Omega)}{\theta[\widetilde{P}_m + \sum_{i \in V} \widetilde{P}_i-\tilde{\Delta}](\Omega)^2}\right).
	\end{align}
	
	Let now 
	\begin{equation}
	T_j = V \cup \{j\},
	\end{equation}
	for $j = l,m$, and replace the notation $\widetilde{P}_i$ with the notation $\eta_i$, using our convention for sums:
	\begin{align}
	[a_l,a_m,a_k,a_{\infty}] 
	& = \left( \frac{\theta[\eta_{T_l}-\eta_k-\tilde{\Delta}](\Omega)\theta[\eta_{T_l}+\eta_k-\tilde{\Delta}](\Omega)}{\theta[\eta_{T_l}-\tilde{\Delta}](\Omega)^2}\right) \\  \notag
	& \quad \div \left(\frac{\theta[\eta_{T_m}-\eta_k-\tilde{\Delta}](\Omega)\theta[\eta_{T_m}+\eta_k-\tilde{\Delta}](\Omega)}{\theta[\eta_{T_m}-\tilde{\Delta}](\Omega)^2}\right) \\
	& = \left( \frac{\theta[\eta_{T_l}+\eta_k-\tilde{\Delta}-2\eta_k](\Omega)\theta[\eta_{T_l}+\eta_k-\tilde{\Delta}](\Omega)}{\theta[\eta_{T_l}-\tilde{\Delta}](\Omega)^2}\right) \\  \notag
	& \quad \div \left(\frac{\theta[\eta_{T_m}+\eta_k-\tilde{\Delta}-2\eta_k](\Omega)\theta[\eta_{T_m}+\eta_k-\tilde{\Delta}](\Omega)}{\theta[\eta_{T_m}-\tilde{\Delta}](\Omega)^2}\right).
	\end{align}
	
	To make our expressions shorter, we write
	\begin{equation}
	S_j = T_j \cup \{k\} = V\cup\{j,k\}
	\end{equation}
	for $j = l,m$, so that we have
	\begin{equation}
	\eta_{T_j}+\eta_k = \eta_{S_j},
	\end{equation}
	since $k \not \in T_j$. Continuing our computation we have
	\begin{align}
	[a_l,a_m,a_k,a_{\infty}] 
	& = \left( \frac{\theta[\eta_{S_l}-\tilde{\Delta}-2\eta_k](\Omega)\theta[\eta_{S_l}-\tilde{\Delta}](\Omega)}{\theta[\eta_{T_l}-\tilde{\Delta}](\Omega)^2}\right)   \div \left(\frac{\theta[\eta_{S_m}-\tilde{\Delta}-2\eta_k](\Omega)\theta[\eta_{S_m}-\tilde{\Delta}](\Omega)}{\theta[\eta_{T_m}-\tilde{\Delta}](\Omega)^2}\right).
	\end{align}
	
	We now notice that for $j = l,m$, the characteristics
	\begin{equation}
	\eta_{S_j}-\tilde{\Delta}-2\eta_k  \quad \text{and} \quad \eta_{S_j}-\tilde{\Delta}
	\end{equation}
	differ by an integer vector, namely $-2\eta_k$. Therefore we may apply the quasi-periodicity property of the Riemann theta constant with characteristic given in equation \eqref{eq:ThetaCharQuasiperiodic} to obtain
	\begin{equation}
	\theta[\eta_{S_j}-\tilde{\Delta}-2\eta_k](\Omega)= \exp(4 \pi i (\tilde{\Delta}-\eta_{S_j})_1(\eta_k)_2)\theta[\eta_{S_j}-\tilde{\Delta}](\Omega).
	\end{equation}
	
	Therefore we have
	\begin{align}
	[a_l,a_m,a_k,a_{\infty}] 
	& = \left( \frac{\exp(4 \pi i (\tilde{\Delta}-\eta_{S_l})_1(\eta_k)_2)\theta[\eta_{S_l}-\tilde{\Delta}](\Omega)^2}{\theta[\eta_{T_l}-\tilde{\Delta}](\Omega)^2}\right) \\  \notag
	& \quad \div \left(\frac{\exp(4 \pi i (\tilde{\Delta}-\eta_{S_m})_1(\eta_k)_2)\theta[\eta_{S_m}-\tilde{\Delta}](\Omega)^2}{\theta[\eta_{T_m}-\tilde{\Delta}](\Omega)^2}\right)\\ \notag
	& = \frac{\exp(4 \pi i (\tilde{\Delta}-\eta_{S_l})_1(\eta_k)_2)}{\exp(4 \pi i (\tilde{\Delta}-\eta_{S_m})_1(\eta_k)_2)}\left( \frac{\theta[\eta_{S_l}-\tilde{\Delta}](\Omega)\theta[\eta_{T_m}-\tilde{\Delta}](\Omega)}{\theta[\eta_{S_m}-\tilde{\Delta}](\Omega)\theta[\eta_{T_l}-\tilde{\Delta}](\Omega)}\right)^2.
	\end{align}
	
	We can simplify the sign:
	\begin{align}
	\frac{\exp(4 \pi i (\tilde{\Delta}-\eta_{S_l})_1(\eta_k)_2)}{\exp(4 \pi i (\tilde{\Delta}-\eta_{S_m})_1(\eta_k)_2)} & = \exp(4 \pi i (\eta_{S_m}-\eta_{S_l})_1(\eta_k)_2)\\ \notag
	& = \exp(4 \pi i (\eta_{m}-\eta_{l})_1(\eta_k)_2).
	\end{align}
	
	We now handle the quantity $\tilde{\Delta}$. First, we note that since $\tilde{\Delta}$ is a vector with half-integer entries, $\tilde{\Delta}$ and $-\tilde{\Delta}$ differ by a vector with integer entries. Furthermore, as noted in equation \eqref{eq:Uset}, $\eta_{U_{\eta}}$ and $\tilde{\Delta}$ differ by a vector with integer entries. Therefore $-\tilde{\Delta}$ and $\eta_{U_{\eta}}$ differ by a vector with integer entries, say $n$:
	\begin{equation}
	-\tilde{\Delta} = \eta_{U_{\eta}} + n.
	\end{equation}
	
	Recalling our notation for the symmetric difference of two sets given in equation \eqref{eq:circ}, we have that
	\begin{equation}
	\eta_{S_j} - \tilde{\Delta} = \eta_{S_j}+ \eta_{U_\eta} +n = \eta_{S_j \circ U_\eta} + 2 \eta_{S_j \cap U_\eta} + n,
	\end{equation}
	and
	\begin{equation}
	\eta_{T_j} - \tilde{\Delta} = \eta_{T_j} + \eta_{U_\eta} +n= \eta_{T_j \circ U_\eta} + 2 \eta_{T_j \cap U_\eta}+n,
	\end{equation}
	for $j = l,m$.

	Carrying on with our computation we therefore have
	\begin{align}
	[a_l,a_m,a_k,a_{\infty}]= \exp(4 \pi i (\eta_{m}-\eta_{l})_1(\eta_k)_2)\left( \frac{\theta[\eta_{S_l \circ U_\eta}+2\eta_{S_l \cap U_\eta}+n](\Omega)\theta[\eta_{T_m \circ U_\eta}+2\eta_{T_m \cap U_\eta}+n](\Omega)}{\theta[\eta_{S_m \circ U_\eta}+2\eta_{S_m \cap U_\eta}+n](\Omega)\theta[\eta_{T_l \circ U_\eta}+2\eta_{T_l \cap U_\eta}+n](\Omega)}\right)^2.
	\end{align}
	
	Lastly, we apply the quasi-periodicity property of the Riemann theta constant with characteristic once more to remove the integer vectors appearing in each characteristic. This time around, we note that since all of characteristics appearing above are half-integers, the sign $\exp(2\pi i x_1m_2)$ from the transformation formula will be $\pm 1$. Since all of the theta constants are now squared in the formula, the signs vanish and we finally obtain:
	\begin{equation}
	[a_l,a_m,a_k,a_{\infty}]= \exp(4 \pi i (\eta_{m}-\eta_{l})_1(\eta_k)_2)\left( \frac{\theta[\eta_{S_l \circ U_\eta}](\Omega)\theta[\eta_{T_m \circ U_\eta}](\Omega)}{\theta[\eta_{S_m \circ U_\eta}](\Omega)\theta[\eta_{T_l \circ U_\eta}](\Omega)}\right)^2.
	\end{equation}
	
	This completes the proof
\end{proof}

To finish the proof of Proposition \ref{prop:MainProp}, it remains now only to rewrite it so that the characteristics agree with Takase's and to verify that the signs agree. Indeed, the cross-ratio we compute here in this article agrees with the quotient computed by Takase, since in his article, Takase assumes that $a_{\infty} = \infty$. In that case, we have that
\begin{equation}
[a_l,a_m,a_k,a_{\infty}] = \frac{a_k-a_l}{a_k-a_m}.
\end{equation}

We therefore turn our attention to the characteristics: Following Takase's notation, let $W$ be the complement of $V \cup \{k,l, m,\infty\}$ in $\{1,2,\ldots, 2g+2\}$. Then from the definitions it follows that
\begin{equation}
S_j = V \cup \{k,j\},
\end{equation}
for $j = l,m$. We also have that $T_l\cup\{\infty\}$ is the complement of $W \cup \{k,m\}$ in $\{1,2,\ldots, 2g+2\}$, and $T_m\cup\{\infty\}$ is the complement of $W \cup \{k,l\}$. As a result,
\begin{equation}
((T_m\cup\{\infty\}) \circ U_{\eta})^c = U_{\eta} \circ (W \cup \{k,l\}),
\end{equation}
and
\begin{equation}
((T_l\cup\{\infty\}) \circ U_{\eta})^c = U_{\eta} \circ (W \cup \{k,m\}).
\end{equation}

Now by definition, we have that
\begin{equation}
\eta_\infty = 0,
\end{equation}
since $P_\infty$ is chosen to be the base point of the Abel-Jacobi map, so the divisor $P_\infty-P_\infty$ maps to the identity in $J(C)$. Therefore we have
\begin{equation}
\eta_{(T_j \cup \{\infty\})\circ U_{\eta}} = \eta_{T_j\circ U_\eta},
\end{equation}
for $j = l,m$, since the difference between the two sides of the equality is $\eta_\infty$ which is zero. By equation \eqref{eq:complements}, we have that
\begin{equation}
\eta_{(T_l \cup \{\infty\})\circ U_{\eta}} = \eta_{U_{\eta} \circ (W \cup \{k,m\})}
\end{equation}
and 
\begin{equation}
\eta_{(T_m \cup \{\infty\})\circ U_{\eta}} = \eta_{U_{\eta} \circ (W \cup \{k,l\})}.
\end{equation}

Putting all of this together, we obtain
\begin{equation}
[a_l,a_m,a_k,a_{\infty}]= \exp(4 \pi i (\eta_{m}-\eta_{l})_1(\eta_k)_2)\left( \frac{\theta[\eta_{U_\eta\circ(V \cup \{k,l\})}](\Omega)\theta[\eta_{U_{\eta} \circ (W \cup \{k,l\})}](\Omega)}{\theta[\eta_{U_\eta \circ (V \cup \{k,m\})}](\Omega)\theta[\eta_{U_{\eta} \circ (W \cup \{k,m\})}](\Omega)}\right)^2.
\end{equation}

To verify that the signs agree, we note that before simplifying his expression, Takase has the sign written as 
\begin{equation}
(-1)^{4 (\eta_k)_1 (\eta_l+\eta_m)_2} = \exp(4 \pi i (\eta_k)_1 (\eta_l+\eta_m)_2).
\end{equation}
We first begin by noting that the sign that we obtain is equal to
\begin{equation}
\exp(4 \pi i (\eta_{m}-\eta_{l})_1(\eta_k)_2) = \exp(4 \pi i (\eta_{l}+\eta_{m})_1(\eta_k)_2)
\end{equation}
(here we have replaced the difference of $\eta_m$ and $\eta_l$ with its sum), since both $\eta_l$ and $\eta_m$ have half-integer entries; the sign of the expression depends only on the number of half-integer (as opposed to integer) coordinates, not on their sign or size. Define
\begin{equation}
e_2(\xi,\zeta) = \exp(4 \pi i(\xi_1\zeta_2-\xi_2\zeta_1)),
\end{equation}
then we have
\begin{align}
\exp(4 \pi i (\eta_{l}+\eta_{m})_1(\eta_k)_2)&\exp(4 \pi i (\eta_k)_1 (\eta_l+\eta_m)_2)\\ \notag
& =  \exp(4 \pi i (\eta_{l}+\eta_{m})_1(\eta_k)_2)\exp(-4 \pi i (\eta_k)_1 (\eta_l+\eta_m)_2) \\ \notag
&= e_2(\eta_l+\eta_m,\eta_k) \\ \notag
& = e_2(\eta_l,\eta_k)e_2(\eta_m,\eta_k).
\end{align}
Now it is a property of the eta-maps that $e_2(\eta_i,\eta_j) = -1$ whenever $i\neq j$ (see \cite[Lemma 1.4.13]{poor} or \cite[Proposition 3.5]{BILV}), so the expression is $1$. This shows that the two signs (Takase's and ours) are always the same. This completes the proof of Proposition \ref{prop:MainProp}.

We now end with the proof of Corollary \ref{cor:Rosenhain}:
\begin{proof}[Proof of Corollary \ref{cor:Rosenhain}]
	To obtain the values $\lambda_i$, we post-compose the degree $2$ morphism $x \colon C \to \mathbb{P}^1$ with a linear fractional transformation of $\mathbb{P}^1$ sending $x(P_1)$ to $0$, $x(P_2)$ to $1$ and $x(P_{2g+2})$ to $\infty$. This new map is again a degree $2$ morphism $C \to \mathbb{P}^1$, and so the result of Proposition \ref{prop:MainProp} applies. In addition, we use that for this particular map, if $\lambda_j = x(P_j)$, then we have 
	\begin{equation}
	\lambda_j = \frac{0-\lambda_j}{0-1}= [\lambda_j,1,0,\infty] =\frac{x(P_1) - x(P_j)}{x(P_1)-x(P_2)}.
	\end{equation}
	Therefore we fix $k = 1$, $l=2$ and $m= j$ to obtain
	\begin{equation}
	\lambda_j =  \exp(4 \pi i (\eta_j-\eta_2)_1(\eta_1)_2)\left( \frac{\theta[\eta_{U_\eta\circ(V \cup \{1,2\})}](\Omega)\theta[\eta_{U_{\eta} \circ (W \cup \{1,2\})}](\Omega)}{\theta[\eta_{U_\eta \circ (V \cup \{1,j\})}](\Omega)\theta[\eta_{U_{\eta} \circ (W \cup \{1,j\})}](\Omega)}\right)^2.
	\end{equation}
\end{proof}

\bibliographystyle{unsrt}

\end{document}